\documentclass[a4paper,11pt]{article}

\usepackage{xypic,amsmath,amssymb,latexsym,enumerate}
\usepackage{theorem,eepic} 
 \usepackage{geometry,xspace}
\usepackage[curve]{xy}
\usepackage{float}
\usepackage{epsfig}
\usepackage{epic}
\usepackage{overpic}

\newcommand{\cat}{\mathsf{Cat}}

\def\2cat{\mathrm{2-}\cat}

\def\FTop{\mathsf{FTop}}

\newtheorem{example}{Example}[section]
{\theorembodyfont{\rmfamily}}
{\theorembodyfont{\rmfamily}}
{\theorembodyfont{\rmfamily}}
{\theorembodyfont{\rmfamily}}

\newtheorem{thm}[example]{Theorem}

\newtheorem{blank}[example]{\hspace{-0.3em}}
{\theorembodyfont{\rmfamily}}

\newcommand{\sqdiagram}[8]{\vcenter{ \xymatrix { #1 \ar[r]^-{#2} \ar[d]_{#4} &
#3 \ar[d]^{#5} \\ #6 \ar[r]_-{#7} & #8 }}}

\newcommand{\directs}[2]{\def\objectstyle{\scriptstyle} \objectmargin={0pt}
\xy
(0,4)*+{}="a",(0,-2)*+{\rule{0em}{1.5ex}#2}="b",(7,4)*+{\;#1}="c"
\ar@{->} "a";"b" \ar @{->}"a";"c" \endxy }

\newcommand{\xdirects}[2]{\def\objectstyle{\scriptstyle} \objectmargin={0pt}
\xy
(0,0)*+{}="a",(0,-6)*+{\rule{0em}{1.5ex}#2}="b",(7,0)*+{\;#1}="c"
\ar@{->} "a";"b" \ar @{->}"a";"c" \endxy }

\newcommand{\sdirects}[2]{\def\objectstyle{\scriptstyle} \objectmargin={0pt}
\xy
(0,2.2)*+{}="a",(0,-2.5)*+{\rule{0em}{1.5ex}#2}="b",(7,2.2)*+{\;#1}="c"
\ar@{->} "a";"b" \ar @{->}"a";"c" \endxy }

\newcommand{\bl}{\mbox{\rule{0.08em}{1.7ex}\hspace{-0.00em}\rule{0.7em}{0.2ex}}}
\newcommand{\br}{\mbox{\rule{0.7em}{0.2ex}\hspace{-0.04em}\rule{0.08em}{1.7ex}}}
\newcommand{\tr}{\mbox{\rule[1.5ex]{0.7em}{0.2ex}\hspace{-0.03em}\rule{0.08em}{1.7ex}}}
\newcommand{\tl}{\mbox{\rule{0.08em}{1.7ex}\rule[1.54ex]{0.7em}{0.2ex}}}

\newcommand{\hh}{\mbox{\rule{0.7em}{0.2ex}\hspace{-0.7em}\rule[1.5ex]{0.70em}{0.2ex}}}
\newcommand{\vv}{\mbox{\rule{0.08em}{1.7ex}\hspace{0.6em}\rule{0.08em}{1.7ex}}}

\def\epsilon{\varepsilon}

\renewcommand{\ge}{\geqslant}

\def\<{\langle \! \langle}
\def\>{\rangle \! \rangle}

\def\A{\mathsf{A}}

\newcommand{\I}{{\mathcal I}}

\newcommand{\C}{\mathcal{C}}

\def\B{\mathcal{B}}

\def\dom{\operatorname{dom}}
\def\Hol{\operatorname{Hol}}
\def\A{\alpha}
\def\B{\beta}

\def\C{\mathsf{C}}
\def\subset{\subseteq}

\begin{document}
\title{Three
themes in the work of  Charles Ehresmann: \\ Local-to-global;
Groupoids; Higher dimensions.  \thanks{Expansion of an invited
talk given to the 7th Conference on the  Geometry and Topology of
Manifolds:  The Mathematical Legacy of Charles Ehresmann,
 Bedlewo 8.05.2005-15.05.2005 (Poland)}}

\author{Ronald Brown\thanks{Partially supported by
a Leverhulme Emeritus Fellowship 2002-2004}\\ University of Wales,
Bangor.}

\maketitle
\begin{center}
  UWB Math Preprint 06.03
\end{center}
\begin{abstract}
This paper illustrates the themes of the title in terms of: van
Kampen type theorems for the fundamental groupoid; holonomy and
monodromy groupoids; and higher homotopy groupoids. Interaction
with work of the writer is explored. \footnote{
MSClass2000: 01A60,53C29,81Q70,22A22,55P15 }
\end{abstract}

\section{Introduction}

It is a pleasure  to honour Charles Ehresmann by giving  a
personal account of  some of the major themes in his work which
interact with mine. I hope it will be useful to suggest how these
themes are related, how the pursuit of them gave a distinctive
character to his aims and his work, and how they influenced my own
work, through his writings and through other people.

Ehresmmann's work is so extensive that a full review would be  a
great task, which to a considerable extent is covered  by Andr\'ee
Ehresmann in her commentaries in the collected works \cite{Eh-Oe}.
His wide vision is shown by his description of his  overriding aim
as: `To find the structure of everything'.  `To find structure' is
related to the Bourbaki experience and aim, in which he was a
partner. A description of a new structure is in some sense a
development of part of a new language: the aim of doing this
contrasts with that of many, who feel that the development of
mathematics is mainly guided by the  solution of famous problems.

The notion of structure is also related to the notion of
\emph{analogy}.  It one of the triumphs of category theory in the
20th century to make progress in unifying mathematics through the
finding of analogies between the behaviours of structures across
different areas of mathematics.

This theme is elaborated in the article \cite{Br-Po-analogy}. That
article argues that many analogies in mathematics, and in many
other areas, are not between objects themselves but between the
relations between objects. Here we mention only the notion of
\emph{pushout}, which we use later in discussing the van Kampen
Theorem. A pushout has the same definition in different categories
even though the construction of pushouts in these categories may
be wildly different. Thus a concentration on the constructions
rather than on the universal properties may lead to a failure to
see the analogies.

Ehresmann developed new concepts and new language which have been
very influential in mathematics: I mention only fibre bundles,
foliations, holonomy groupoid, germs, jets, Lie groupoids. There
are other concepts whose time perhaps is just coming or has yet to
come: included here might be ordered groupoids, multiple
categories.

In this direction of developing language, we can usefully quote
G.-C. Rota \cite[p.48]{Rota}:
\begin{quote}
  ``What can you prove with exterior algebra that you cannot prove
without  it?" Whenever you hear this question raised about some
new piece of  mathematics, be assured that you are likely to be in
the presence of something important. In my time, I have heard it
repeated for random  variables, Laurent Schwartz' theory of
distributions, ideles and Grothendieck's schemes, to mention only
a few. A proper retort might be: ``You are right. There is nothing
in yesterday's mathematics that could not also be proved without
it. Exterior algebra is not meant to prove old facts, it is meant
to disclose a new world. Disclosing new worlds is as worthwhile a
mathematical enterprise as proving old conjectures.''
\end{quote}
\section{Local-to-global questions}
\label{loctoglob} Ehresmann developed many new themes in category
theory. One example is structured categories, with  principal
examples those of  differentiable categories and groupoids, and of
multiple categories. His work on these is quite disparate from the
general development of category theory in the 20th century, and it
is interesting to search for reasons for this. One must be the
fact that he used his own language and notation. Another is surely
that his early training and motivation came from analysis, rather
than from algebra, in contrast to the origins of category theory
in the work of Eilenberg, Mac~Lane and of course Steenrod, centred
on homology and algebraic topology. Part of the developing
language of category theory became essential in those areas, but
other parts, such as that of algebraic theories,  groupoids,
multiple categories, were not used till fairly recently.

It seems likely that Ehresmann's experience in analysis led him to
the major theme of \emph{local-to-global} questions. I first
learned of this term from Dick Swan in Oxford in 1957-58, when as
a research student I was writing up notes of his Lectures on the
Theory of Sheaves \cite{Swan}. Dick  explained to me that two
important methods for local-to-global problems were sheaves and
spectral sequences---he was thinking of Poincar\'e duality, which
is discussed in the lecture notes,  and the more complicated
Dolbeaut's theorem for complex manifolds. But in truth such
problems are central in mathematics, science and technology. They
are fundamental to differential equations and dynamical systems,
for example. Even deducing consequences of a set of rules is a
local-to-global problem: the rules are applied locally, but we are
interested in the global consequences.

My own work on local-to-global problems arose from writing an
account of the Seifert-van Kampen theorem on the fundamental
group. This theorem can be given as follows, as first shown by
R.H. Crowell:
\begin{thm}{\em \cite{Cr59}}\label{vktgp} Let the space $X$ be the
union of open sets $U,V$ with intersection $W$, and suppose
$W,U,V$ are path connected. Let $x_0 \in W$. Then   the diagram of
fundamental group morphisms induced by inclusions
\begin{equation} \label{push}
{\sqdiagram{\pi_1(W,x_0)}{i}{\pi_1(U,x_0)}{j}{}{\pi_1(V,x_0)}{}{\pi_1(X,x_0)}}
\end{equation}
is a pushout of groups. \end{thm} Here the `local parts' are of
course $U,V$ put together with intersection $W$ and the result
describes completely, under the open set and connectivity
conditions, the nonabelian fundamental group of the global space
$X$. This theorem is usually seen as a necessary part of basic
algebraic topology, but one without higher dimensional analogues.

In writing the first 1968 edition of the book \cite{Br-newbook}, I
noted  that to compute the fundamental group of the circle one had
to develop something of  covering space theory. Although  that is
an excellent subject in its own right, I became irritated by this
detour. After some time, I found  work of Higgins on groupoids,
\cite{Hi-present}, which defined free products with amalgamation
of groupoids, and this led to a more general formulation of
theorem \ref{vktgp} as follows:
\begin{thm}{\em \cite{Br-firstvkt}}\label{vktgpd} Let the space $X$ be the
union of open sets $U,V$ with intersection $W$. Let $X_0$ be a
subset of $W$ meeting each path component of $W$.  Then the
diagram of fundamental groupoid morphisms induced by inclusions
\begin{equation} \label{pushgpd}
{\sqdiagram{\pi_1(W,X_0)}{i}{\pi_1(U,X_0)}{j}{}{\pi_1(V,X_0)}{}{\pi_1(X,X_0)}}
\end{equation}
is a pushout of groupoids. \end{thm} Here $\pi_1(X,X_0)$ is the
fundamental groupoid of $X$ on a set $X_0$ of base points: so it
consists of homotopy classes rel end points of paths in $X$
joining points of $X_0 \cap X$.

In the case $X$ is the circle $S^1$, one chooses $U,V$ to be
slightly extended semicircles including $X_0=\{+1,-1\}$. The point
is that in this case $W=U \cap V$ is not path connected and so it
is not clear where to choose a single base point. The day is saved
by hedging one's bets, and using two base points.

The proof of theorem \ref{vktgpd} uses the same tricks as  to
prove theorem \ref{vktgp}, but in a broader context. In order to
compute fundamental groups from this theorem, one can set up some
general combinatorial groupoid theory, see
\cite{Br-newbook,Hi-book}. A key feature of this theory is the
groupoid $\I$, the indiscrete groupoid on two objects $0,1$, which
acts as a unit interval object in the category of groupoids. It
also plays a r\^{o}le analogous to that of the infinite cyclic
group $\C$ in the category of  groups. One then compares the
pushout diagrams, the first in spaces, the second in groupoids:
\begin{equation}\label{circle}
\vcenter{\xymatrix{\{0,1\} \ar [r] \ar [d] & \{0\} \ar [d]\\
[0,1] \ar [r] & S^1_{\phantom{1}}}}\hspace{5em} \vcenter{
\xymatrix{\{0,1\} \ar [r] \ar [d] & \{0\} \ar [d]\\ \I \ar [r]_u &
\C}}\end{equation} to see how this version of the van Kampen
Theorem gives an analogy between the geometry, and the algebra
provided by the notion of groupoid. This kind of result is seen as
`change of base' in \cite{Br-base-ch}.

The fundamental group is a kind of anomaly in algebraic topology
because of its nonabelian nature. Topologists in the early part of
the 20th century were aware that: \begin{itemize}
 \item the non commutativity of the  fundamental group was useful in
applications; \item for path connected $X$ there was an
isomorphism
$$H_1(X) \cong \pi_1(X,x)^{\mathrm{ab}};   $$
\item the abelian homology groups existed in all dimensions.
\end{itemize}
Consequently there was a desire to generalise the nonabelian
fundamental group to all dimensions.

In 1932 \v{C}ech submitted a paper on higher homotopy groups
$\pi_n(X,x)$ to the ICM at Zurich, but it was quickly proved that
these groups were abelian for $n \ge 2$, and on these grounds
\v{C}ech was persuaded to withdraw his paper, so that only a small
paragraph appeared in the Proceedings \cite{Ce}. We now see the
reason for the commutativity as the result (Eckmann-Hilton) that a
group internal to  the category of groups is just an abelian
group. Thus the vision of a non commutative higher dimensional
version of the fundamental group has since 1932 been generally
considered to be a mirage.

Theorem \ref{vktgpd} is also anomalous: it is a colimit type
theorem, and so yields \emph{complete information} on the
fundamental groups which are contained in it,  even in the non
connected case, whereas the usual method in algebraic topology is
to relate different dimensions by exact sequences or even spectral
sequences, which usually yield information only up to extension.
Thus exact sequences by themselves cannot show that a group is
given as an HNN-extension: however such a description may be
obtained from a pushout of groupoids, generalising the pushout of
groupoids in diagram \ref{circle}.

It was then found that the theory of covering spaces could be
given a nice exposition using the notion of covering morphism of
groupoids. Even later, it was found by Higgins and Taylor
\cite{HiTa} that there was a nice theory of orbit groupoids which
gave models of orbit spaces.

The objects of a groupoid add a `spatial component' to group
theory, which is essential in many applications. This is evident
in many parts of Ehresmann's work. Another view of  this anomalous
success of groupoids is that  they have structure in two
dimensions,  0 with the objects and 1 with the arrows. We have a
colimit type theorem for this larger structure, and so a good
model of the geometry. Useful information on fundamental groups is
carried by the fundamental groupoid.

It is therefore natural to seek  for higher homotopy theory
algebraic models which:\begin{itemize} \item have structure in a
range of dimensions; \item contain useful information on classical
invariants, and \item satisfy  van Kampen type theorems.
\end{itemize}That is, we seek \emph{nonabelian methods for
higher dimensional local-to-global problems in homotopy theory}.
We return to this theme in section \ref{higherdimension}, which
gives an indication of some of the motivation for the writer for
{\em higher dimensional algebra}.

\section{Holonomy and monodromy}
\label{holmon} Once I had been led to groupoids by Philip Higgins,
and having been told by G.W. Mackey in 1967 of his work on ergodic
groupoids, \cite{Mackey}, it was natural to consider topological
groupoids and  differentiable groupoids, and to seek their
properties and applications, \cite{Br-Hardy1,Br-Hardy2}.

So I came across papers of Ehresmann, \cite{Eh}, and of Jean
Pradines, \cite{Pr},  and in 1981 I visited Jean  in Toulouse
under British Council support, to try and understand something of
his papers. I saw the grand vision of the whole scheme of
generalising the relation between Lie groups and Lie algebras to a
relation between Lie groupoids and Lie algebroids, which has now
become a large theory, see for example \cite{Ma}. We concentrated
on  his first paper \cite{Pr}, which states theorems but gives no
indication of proofs. What I found remarkable was first of all the
beautiful constructions Jean explained, which seemed to me clear
in principle, and then the fact that he gave for holonomy and
monodromy universal properties: these are fairly rare in
differential topology.

Jean was interested in the \emph{monodromy principle}, which
involves the following situation:
$$
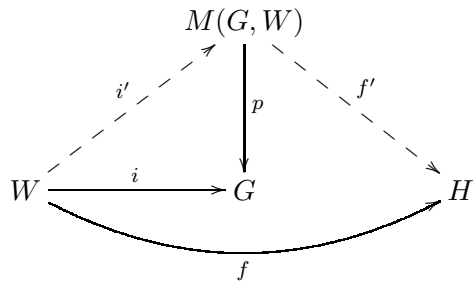
\xymatrix@=4pc{& M(G,W) \ar [d]^p \ar @{-->} [dr] ^{f'}& \\
W \ar @{-->} [ur] ^{i'} \ar [r] ^i \ar @/_2pc/ [rr]_f & G & H }$$
Here $G$ is a topological groupoid, $W$ is an open subset of $G$
which contains the set of identities $G_0$ of $G$, and $i$ is the
inclusion. The aim is to find a topological groupoid $M(G)$,
called the \emph{monodromy  groupoid of $G$}, with morphism of
topological groupoids $p: M(G,W) \to G$ and a `local morphism'
$i': W \to M(G,W)$ such that for any local morphism $f:W \to H$ to
a topological groupoid $H$, $f$ `locally extends' to $f': M(G,W)
\to H$. The existence such $M(G,W), p, i'$ with this universal
property is called the \emph{monodromy principle}.

This idea was expressed in Chevalley's famous book on `Lie
groups', \cite{Ch}, for the case the groupoid $G$ is the trivial
groupoid $X \times X$ on a manifold $X$, when $M(G,W)$ becomes,
for sufficiently small $W$,  the fundamental groupoid of $X$ and
$p: M(G,W) \to G$ is the `anchor map' as defined in Mackenzie's
book \cite{Ma}.  Pradines is the first, I believe,  to see how to
extend it this notion to the  case of a general groupoid, as
announced in \cite{Pr}.

We have to explain the term `local morphism'. Note that if $u,v
\in W$ and $uv$ is defined in $G$, this does not necessarily mean
that $uv \in W$. So the algebraic structure that $W$ inherits from
$G$ is what is called a `local groupoid structure'. A local
morphism is therefore one that preserves the local groupoid
structure.

It is easy to construct  algebraically the groupoid $M(G,W)$ and
$p:M(G,W) \to G$ so that it has the required universal property,
algebraically, as follows. The source, target and identity maps
for $G$ induce on $W$ the structure of a reflexive graph. So one
forms $F(W)$, the free groupoid on this reflexive graph, the
reflexive condition ensuring that identities in $W$ become
identities in $F(W)$. There is a function $j: W \to F(X)$ which
sends $w \in W$ to the corresponding generator $[w] \in F(W)$ and
then one quotients the groupoid $F(W)$ by the relations
$[u][v]=[uv]$ for all $u,v \in W$ such that $uv$ is defined in
$W$. This defines $M(G,W)$, and the function $j$ induces the local
morphism $i'$, which is injective since it is determined by the
inclusion $i$.

The problem is how to topologise $M(G,W)$ so that it becomes a
topological groupoid for which  the monodromy principle is
satisfied not only  algebraically but also with regard to
continuity, or differentiability.

Pradines solved this by a beautiful holonomy construction, which
he explained to me in 1981 during my visit, and which develops
Ehresmann's ideas.

He took the view that the pair $(M(G,W),W')$, where $W'=i'(W)$,
should be regarded as a `locally Lie groupoid', in that the Lie
structure resides only on $W'$. This raised the general problem of
when a locally Lie groupoid is \emph{extendible}, i.e. is obtained
from an open subset containing the identities of a Lie groupoid.
(The term used in \cite{Pr} is {\em morceau d'un groupo\"{\i}de
diff\'erentiable}.)

We therefore start again and consider a groupoid $G$ and a subset
$W$ of $G$ such that $W$ contains the identities of $G$ and $W$
has the structure of topological space and even of a manifold. We
ask: what conditions should be put on $W$ which ensure that $G$
can be given the structure of topological or Lie groupoid for
which $W$ is an open set? In short, we ask is the pair $(G,W)$
extendible?

This is a classical question in the case $G$ is a group, and the
answer is given in books on topological groups. A topology on $G$
is obtained by taking as subbase the sets $gU$ for all $U$ open in
$W$ and $g \in G$. The conditions for this to give a topological
group  structure on $G$ are fairly mild, and are roughly  that the
algebraic operations on $G$ should be as continuous on $W$ as can
be expected given that these operations are defined only partially
on $W$.

One of the reasons this works is that in a topological group, the
left multiplication operator $L_g: G \to G$,  given for $g \in G$
by $u \mapsto gu$,  is a homeomorphism and so maps open sets to
open sets.

This property of  $L_g$  no longer holds  when $G$ is a
topological groupoid, because $L_g$ is not defined on all of $G$.
This reflects the considerable change in moving from groups to
groupoids, that is, from algebraic operations  always defined to
those only  partially defined. This is not a loss of information:
the wider concept has greater powers of expression.

There is also a wide range of examples where $(G,W)$ is not
extendible, many coming from the theory of foliations. For
example, if $\cal F$ is the foliation of the M\"{o}bius Band $M$
given by circles going once or twice round the band, and $R$ is
the equivalence relation whose equivalence classes are the leaves,
then $R$ is a subset of $M \times M$ and as such is a topological
groupoid. But $R$ is not a submanifold of $M \times M$, since it
has self intersections. However the foliation structure does
determine a locally Lie groupoid $(R,W)$, although $W$ is not open
in $R$. The argument for  these facts is spelled out in
\cite{BrMu-fol}.

Given a groupoid $\alpha,\beta:G \rightrightarrows G_0$, then  an
\emph{admissible section $s: G_0 \to G$  of $\alpha$} satisfies
$\alpha \circ s=1$ and $\beta \circ s: G_0 \to G_0$ is a
bijection.  We follow Mackenzie in \cite{Ma} in calling $s$ a
\emph{bisection}. We can also regard $s$ as a homotopy in the
category of groupoids $1 \simeq a$ where $a:G \to G$ is an
automorphism of $G$; here homotopy is defined by the `unit
interval object' $\I$. This interpretation has intuitive value,
and is suggestive for analogues in higher dimensions, as applied
in \cite{Br-Ic}.

Now suppose $G$ is as above but $G_0$ is also a topological space.
Then a \emph{local bisection} of $G$ is a function $s:U \to G$
with $U$ open in $G_0$ such that  $\alpha \circ s=1_U$ and $\beta
\circ s$ maps $U$ homeomorphically to its image which is also open
in $G_0$. There is an `Ehresmannian composition' $s*t$ of local
bisections. We first make clear that if $g: x \to y$ and $h: y \to
z$ in $G$ then their composition in $G$ is $hg: x \to z$. So the
composition of local bisections is $$(s*t)(x)= (s\beta t(x))
(tx)$$ for $x \in G_0$. This means that in general the domain of
$s*t$ is smaller than that of $t$, and may even be empty. This
composition makes the set of local bisections into an
\emph{inverse semigroup}. Recall that this is a semigroup in which
for each element $s$ there is  a unique element $s'$, called a
{\em relative inverse for $s$},  such that $s'ss'=s', ss's=s$.
Pseudogroups, a concept first defined by Ehresmann in
\cite{Eh-pseudo}, give examples of such structures. We write this
inverse semigroup as $\Gamma(G)$: of course it depends on the
topology of $G_0$.

A left partial `adjoint' operation $L_s$ of the local bisection
$s$ on $G$ is defined by $L_s(g)= (s\beta g) g, g \in G$. It is
easy to prove that if $G$ is a topological groupoid and $s$ is a
continuous local bisection, then $L_s$ maps open sets of $G$ to
open sets of $G$. Thus the important observation is that for
adjoint operations on topological or Lie groupoids it is not
enough to rely on the elements of $G$: we need the local
continuous or smooth bisections, which are kind of `tubes' rather
than `elements',  to transport the local structure of $G$.

Now suppose given a pair $(G,W)$ such that  $G$ is a groupoid,
$G_0 \subset W \subset G$, $W$ is a manifold, and the groupoid
operations are `as smooth as possible' on $W$. By a \emph{smooth
local  bisection of $(G,W)$} we mean a local bisection $s$ of $G$
such that $s$ takes values in $W$ and is smooth. The set of smooth
local bisections forms a subset $\Gamma^{(r)}(W)$ where $r$
denotes the class of differentiability of the manifold $M=G_0$.

Pradines' key definition is to form the sub-inverse semigroup
$\Gamma^{(r)}(G,W)$ of $\Gamma(G)$ \emph{generated by}
$\Gamma^{(r)}(W)$. My interpretation  is that an element of
$\Gamma^{(r)}(W)$ can be thought of as a \emph{local procedure}
and an element of $\Gamma^{(r)}(G,W)$ can be thought of as an
\emph{iteration of local procedures}. Thus an iteration of local
procedures need not be local, and this is one of the basic
intuitions of non trivial holonomy.

We say that $(G,W)$ is \emph{sectionable} if for all $w \in W$
there is a smooth local bisection $s$ in $W$ whose domain includes
$\alpha w$ and with $s\alpha w=w$.

The next step is to form from $\Gamma^{(r)}(G,W)$ the associated
sheaf of germs $J^{(r)}(G,W)$: the elements of this are written
$[s]_x$ where $s \in \Gamma^{(r)}(G,W)$ and $x \in \dom s$. The
inverse semigroup structure on  $\Gamma^{(r)}(G,W)$ induces a
groupoid structure on $J^{(r)}(G,W)$. This contains a subgroupoid
$J^0$ whose elements are germs $[s]_x$ such that $\beta sx= x$ and
there is a neighbourhood $U$ of $x$ such that $s|U \in
\Gamma^{(r)}(W)$; in words, $s$ is an iteration of local
procedures about $x$ which is still a local procedure. It is a
proposition that $J^0(G,W)$ is a normal subgroupoid of
$J^{(r)}(G,W)$. The {\it holonomy groupoid} $\Hol^{(r)}(G,W)$ is
defined to be the quotient groupoid $J^{(r)}(G,W)/ J^0$. The class
of $[s]_x$ in the holonomy groupoid is written $\langle s
\rangle_x$. There is a projection $p: \Hol(G,W) \to G$ given by $
\langle s \rangle_x \mapsto s(x)$.

 The intuition is that first of all   $W$    embeds  in $\Hol^{(r)}(G,W)$, by
$w \mapsto \langle f \rangle_{\alpha w}$, where $f$ is a local
smooth bisection such that $f \alpha w=w$,  and second that
$\Hol(G,W)$ has enough local sections for it to obtain a topology
by translation of the topology of $W$.

Let  $s\in \Gamma^c  (G,W)$.  We define a partial function
$\sigma_s  :W\to  \Hol^{(r)}(G,W)$. The domain of  $\sigma_s$   is
the set of  $w\in W$ such that  $\beta w\in  \dom  (s)$. The value
$\sigma_s w$  is obtained as follows.  Choose a  smooth bisection
$f$ through $w$. Then we set
 $$\sigma_s w=\langle s\rangle_{\B w}
\langle f\rangle_{\A w} =\langle sf\rangle _{\A w}  .$$ Then
$\sigma_s w$ is independent of the choice of the local section $f
$. It is proved in detail in \cite{Ao-Br} that these $\sigma_s$
form a set of charts for $\Hol^{(r)}(G,W)$ making it into a Lie
groupoid with a universal property.

The books \cite{Ma,MoMr} argue that the most efficient treatment
of the holonomy groupoid of a foliation is via the monodromy
groupoid, which is itself defined using the fundamental groupoid
of the leaves. However there are  counter arguments.

One fact is that their arguments do not so far obtain a monodromy
principle, which is obtained by the opposite route in
\cite{Pr,Br-Mu}. Thus there is loss of a universal principle, with
its potentiality for enabling analogies.

The second problem is that the route through the fundamental
groupoid is based on paths, and so on the standard notion of a
topological space, and its exemplification  as a manifold.

The Pradines' approach gives a clear realisation of the intuitive
idea of \emph{ iteration of local procedures}, without requiring
the notion of path to `carry' these procedures, as happens for
example in the usual process of analytic continuation. It is
possible that this idea would lead to wider applications of non
abelian groupoid like methods for local-to-global problems. For
example, the following picture illustrates a chain of local
procedures from $a$ to $b$:
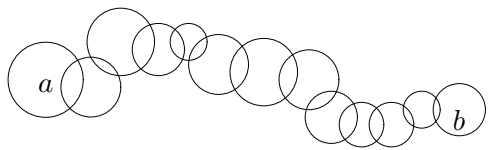
\begin{center}

 \setlength{\unitlength}{10mm} \linethickness{0.4pt}
\begin{picture}(18,-10)(-4,0)
\put(1,0){\makebox(0,0)[t]{$a$}} \put(1,0){\circle{1}}
\put(2,0.5){\circle{0.9}} \put(1.6,-0.1){\circle{0.8}}
\put(2.5,0.4){\circle{0.7}} \put(2.9,0.5){\circle{0.5}}
\put(3.3,0.2){\circle{0.8}} \put(3.9,0.1){\circle{0.9}}
\put(4.5,0){\circle{0.8}} \put(4.8,-0.5){\circle{0.7}}
\put(5.2,-0.6){\circle{0.6}} \put(5.6,-0.6){\circle{0.6}}
\put(6,-0.4){\circle{0.5}} \put(6.5,-0.4){\circle{0.7}}
 \put(6.5,-0.4){\makebox(0,0)[t]{$b$}}

\end{picture}
\end{center}
\vspace{3ex}
 We would like to be able to define such a chain, and
equivalences of such chains, without resource to the notion of
`path'. The claim is that a candidate for this lies in the
constructions of Ehresmann and Pradines for the holonomy groupoid.

Here are some final questions in this area.

Can one use these ideas in other situations to obtain monodromy
(i.e. analogues of `universal covers') in situations where paths
do not exist but `iterations of local procedures' do? It is even
possible that holonomy may exist in wider situations, but not
monodromy.

How widely useful is the notion of `locally Lie groupoid' as a
context for describing local situations? Is there a locally Lie
groupoid $(G,W)$ where $G$ is an action groupoid?

A further point is that Pradines was very keen on having a theory
for germs of locally Lie groupoids. It is in these terms that the
theorems in \cite{Pr} are stated. Such work is not considered in
\cite{Ao-Br,Br-Mu}.

One could consider other structures on $W$ and ask for bisections
which preserve these structures. Simple examples for degrees of
differentiability, due to Pradines \cite{Prp}, are given in
\cite{Ao-Br}. One might consider other geometric structures, such
as  Riemannian, or Poisson.  \.{I}\c{c}en and I, in looking for
notions of double holonomy,  have considered in \cite{Br-Ic} a
double groupoid $G$ and linear bisections.

\section{Higher dimensions} \label{higherdimension}
After writing out the proofs of theorem  \ref{vktgpd} a number of
times to make the exposition clear, it became apparent to me in
1966 that the method of proof ought to extend to higher
dimensions, if one had the right gadgets. So this was an idea of a
proof in search of a theorem. The first search was for a higher
dimensional version of the fundamental groupoid on a set of base
points.

It was at this point that I found Ehresmann's book, \cite{Ehbook}.
It was difficult to understand, but it did contain a definition of
double category, and a key example, the double category $\Box C$
of commuting squares in a category $C$. What caused problems for
me was to find any construction of a fundamental double groupoid
of a space which really contained 2-dimensional homotopy
information. In fact a solution to this problem was not published
till 2002, \cite{Br-Har-Ka-Po}. The construction of higher
homotopy groupoids went in 1975 in a different direction, using
pairs of spaces,  and then filtered spaces, in work with Philip
Higgins, which is surveyed in \cite{Br-fields}.

As an intermediate step, I decided in 1970 to investigate the
notion of double groupoid purely algebraically, to see whether the
putative homotopy double groupoid functor would take values in a
category of some interest.

Work with Chris Spencer, \cite{Br-Sp1}, found a relationship
between double groupoids and the crossed modules introduced by
J.H.C. Whitehead to discuss the second relative homotopy group and
its boundary $\partial: \pi_2(X,A,a) \to \pi_1(A,a)$. We found a
functor from crossed modules to a certain kind of double groupoid,
which was later called {\em edge symmetric}, in that the groupoids
of vertical and of horizontal edges are the same. The next
question was what kind of double groupoids arose in this way?

At the same time we were looking at conjectured proofs of a
vaguely formulated  2-dimensional  van Kampen type theorem. It was
clear that we needed to generalise commutative squares to
commutative cubes in the context of double groupoids, in such a
way that any composite of commutative cubes is commutative.

A square in a groupoid
$$ \def\labelstyle{\textstyle}\vcenter
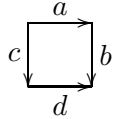
{\xymatrix@M=0pt @=2pc{
\ar  [r]^a \ar  [d]_c &\ar[d]^b\\
\ar [r]_d&}} $$ is commutative if and only if $ab=cd$, or,
equivalently, $a=cdb^{-1}$. We wanted an analogous  expression to
the last in the case of a cube. If we fold flat five of the faces
of a cube we get a net such as the first of the following:
$$\def\labelstyle{\textstyle}\vcenter
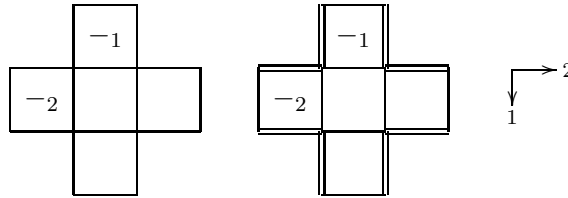
{\xymatrix@M=0pt @=2pc{
 &  \ar @{-}
[ddd] \ar@{-} [r]\ar @{} [dr] | {-_1} & \ar
@{-} [ddd] &  \\
\ar @{} [dr] | {-_2}\ar @{-} [rrr]\ar @{-} [d] &&&\ar@{-}[d]\\
\ar @{-} [rrr]&&&\\
 & \ar@{-} [r]& &}} \qquad \vcenter{\xymatrix@M=0pt @=2pc{
 &  \ar @{=} [d] \ar@{-} [r]\ar @{} [dr] | {-_1} & \ar
@{=} [d] &  \\
\ar @{} [dr] | {-_2}\ar @{=} [r]\ar @{-} [d] & \ar @{-} [d] \ar @{-} [r]&\ar@{-} [d]\ar @{=} [r] &\ar@{-}[d]\\
\ar @{=} [r] &\ar@{-} [r]\ar @{=} [d]&\ar @{=} [d]\ar@{=} [r]&\\
 & \ar@{-} [r]& &}}\qquad   \directs{2}{1} $$
which is not a composable set of squares, since only rectangular
decompositions are composable.  (The third graphic indicates that
squares should have compositions in two directions, and the
symbols $-_1,-_2$ indicate inverses of these compositions.)
However in the second diagram we have noted by double lines
adjacent pairs of edges which are the same. Therefore we must
assume  we have canonical elements, called {\it connections},  to
fill the corner holes as in the diagram:
$$\def\labelstyle{\textstyle}
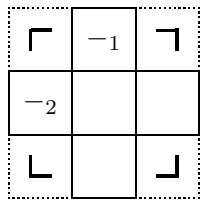
\xymatrix@M=0pt @=2pc{\ar @{..} [r]
\ar@{..} [d] \ar @{} [dr] | \tl &  \ar @{} [dr] | {-_1}\ar @{-}
[ddd] \ar@{-} [r] & \ar
@{-} [ddd] \ar@{..} [r] \ar @{} [dr] | \tr & \ar@{..} [d]\\
 \ar @{} [dr] | {-_2}\ar @{-} [rrr]\ar @{-} [d] &&&\ar@{-}[d]\\
\ar @{-} [rrr]\ar @{..} [d] \ar @{} [dr] | \bl&&\ar @{} [dr] |
\br&\ar@{..} [d]\\
\ar @{..} [r] & \ar@{-} [r]& \ar@{..} [r]&}$$ The rules for
connections include the {transport law:}
$$\def\labelstyle{\textstyle}\vcenter
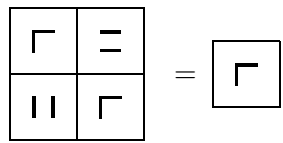
{ \xymatrix@M=0pt @=0.7pc{ \ar @{-} [dddd] \ar @{-} [rrrr]
 &&\ar @{-} [dddd] && \ar @{-}
[dddd]  \\& \tl & & \hh&   \\ \ar @{-} [rrrr] &&&&&
\\ & \vv& & \tl &  \\ \ar @{-} [rrrr]   &&&&}}=\; \vcenter{\xymatrix@M=0pt @=0.7pc{
\ar @{-} [dd]
\ar @{-} [rr]  && \ar @{-} [dd] \\& \tl &\\
\ar @{-} [rr]&& }}
$$
borrowed from a law for path connections found in Virsik's paper
\cite{Vi}. For the full story, see \cite{Br-Sp1}.

Virsik was a student of Ehresmann. He left Czechoslovakia at the
time of the Russian invasion, and went to Australia. There he met
Kirill Mackenzie, and suggested to him Lie groupoids and Lie
algebroids as a  PhD topic. Hearing of Kirill's work, I managed to
get British Council  support for him to visit Bangor in 1986, and
this led to his continuing to work in the UK.

In the paper \cite{Vi},  Virsik introduces the notion of path
connection on a principal bundle $p: E \to B$ as follows. Let $G$
be the Ehresmann groupoid $EE^{-1}\rightrightarrows B$  of the
principal bundle. Thus $G(b,b')$ can be identified with the bundle
maps $E_b \to E_{b'}$ of the fibres over $b,b'$ respectively; $G$
may also be identified with the quotient of $E \times _B E $ by
the action of the group of the bundle. Let $\Lambda B$ be the
category of Moore paths on $B$, and let $\Lambda G$ be the
category of Moore paths on $G$. These may be combined easily to
give the structure of a double category.

Virsik defines a {\em path connection} for the bundle $E$ to be a
function $\Gamma: \Lambda B \to \Lambda G$ satisfying two
conditions. One is invariance under reparametrisation, and the
other is known as the transport law. Chris Spencer and I were
amazed when this law was found to be exactly right as a law on the
connections we needed to define a commutative 3-cube.  So we found
an equivalence between crossed modules and what we then called
{\em special double groupoids with special connections},
\cite{Br-Sp1}. These connections, and their generalisation to all
dimensions,  have been central in work on cubical higher groupoids
ever since.

Ehresmann had earlier shown in \cite{Eh-quintettes} that a
2-category gave rise to a double category of `quintettes'
$(u:a,b,c,d)$ where $a,b,c,d$ form a square and $u:cd \Rightarrow
ab$ is a 2-cell. Spencer showed in \cite{Spencer} that this gave
an equivalence between 2-categories and edge symmetric double
categories with connections. See also \cite{Br-Mo}. This has been
generalised to all dimensions in \cite{Agl-Br-St}, and there has
been recent further work on commutative cubes by Higgins in
\cite{Hi-thin} and Steiner \cite{Steiner}.

What has yet to be accomplished is to find relations between the
higher order connections in differential geometry and the geometry
of cubes.

One of the advantages of cubes for local-to-global questions is
that they have an easy to define notion of multiple composition.
Analogous ideas in the globular case present conceptual and
technical difficulties. Multiple compositions allow easily an {\em
algebraic inverse to subdivision}, and this is a key to certain
local-to-global results, as outlined for example in
\cite{Br-fields}. However the algebraic relations between the
cubical, globular, and (in the groupoid case) crossed complexes,
are an essential part of the picture.

It is possible to dream of a unification of  all these themes in a
way which I believe Ehresmann would have favoured, but there still
seems quite a  way to go to realise this. Work is in progress with
Jim Glazebrook and Tim Porter on aspects of this, \cite{Br-Gl-Po},
and many others are working on stacks, gerbes and higher order
groupoids.

\noindent {\it Acknowledgements}: Thanks to the Conference
organisers for support, and to Jim Stasheff for helpful comments
on this paper.

\end{document}